\documentclass[12pt]{article}
\usepackage{amsxtra,amssymb,amsthm,amsmath,latexsym}

\theoremstyle{plain}

\newtheorem{theorem}{Theorem}[section]

\def\oH{\buildrel\circ\over H}
\def\oH1{\buildrel\circ\over H\kern-.02in{}^1}
\def\qed{{\hfill $\Box$}}

\def\const{\hbox{\,const\,}}

\begin{document}


\title{A characterization of unbounded Fredholm operators.
}

\author{
A.G. Ramm\\
LMA/CNRS, 31 Chemin Joseph Aiguier\\
Marseille 13402, cedex 20\\
France\\
 and Mathematics Department, Kansas State University, \\
 Manhattan, KS 66506-2602, USA\\
ramm@math.ksu.edu\\
}

\date{}

\maketitle\thispagestyle{empty}


\section{Introduction and the result.}
This paper is a continuation of \cite{4}, where bounded Fredholm operators
are studied. The theory of bounded linear Fredholm-type operators is
presented in many texts, see e.g. \cite{1}, \cite{2}. This paper is written 
for a broad audience and the author tries to give simple and short 
arguments.

{\it We call a linear closed densely defined operator $A: X \to Y$ acting 
from a
Banach space $X$ into a Banach space $Y$ a Fredholm operator, and write
$A \in Fred (X, Y)$ if and only if
$$R(A) = \overline{R(A)} \eqno{(1.1)}$$
and
$$n(A) = n(A^\ast) := n < \infty, \quad n(A) := \dim N(A),
  \eqno{(1.2)}$$
where $N(A) := \{u: Au =0, u \in D(A)\}$.}

In the literature the Noether operators are sometimes called Fredholm
operators. The Noether operators are operators for which (1.1) holds,
$n(A)< \infty$, $n(A^\ast) < \infty$, but $n(A)$ may not be equal to
$n(A^\ast)$.
Thus $Fred(X,Y)$ is a proper subset of the Noether operators.

The Noether operators are called in honor of F. Noether, who was the first to
study a class of singular integral equations with operators of this class in
1921 \cite{3}.

In \cite{4} a simple and short proof of the Fredholm alternative and
a characterization of Fredholm operators are given for
bounded linear operators. Recall that a linear bounded operator $F$ is called
a finite-rank operator if $\dim R(F) < \infty$, where $R(F)$ is the range
of $F$.

In the present paper the results of \cite{4} are generalized to
the case of closed unbounded linear operators. Namely, the following result
is proved:

\begin{theorem} 
If $A$ is a Fredholm operator 
then
$$A = B-F, \eqno{(1.3)}$$
where $B$ is a linear closed operator, $D(B) = D(A)$, $R(B) = Y$,
$N(B) = \{0\}$, and $F$ is a finite-rank operator.
Conversely, if (1.3) holds, where $B:X\to Y$ is a linear closed densely defined
operator, $R(B) = Y$, $N(B) = \{0 \}$, and $F$ is a finite-rank operator,
then $A$ is closed, $D(A) = D(B)$, and (1.1) and (1.2) hold, so $A$ is a 
Fredholm operator.
\end{theorem}

In section 2 a proof of Theorem 1 is given. In the literature the case of
unbounded Fredholm operators is usually not discussed directly. In \cite{5}
singularities of the parameter-dependent Fredholm operators are studied, and
in \cite{6} applications of the Fredholm operators in branching theory are
presented. Theorem 1.1 is useful, for example, in the theory of 
elliptic boundary value problems, but we do not go into further 
detail (see, e.g.,  \cite{1}, \cite {2}, \cite {6}).

\section{Proof}
1. Assume that $A : X\to Y$ is linear, closed, densely defined
operator, and (1.1) and
(1.2) hold. Let us prove that then (1.3) holds, $D(B) = D(A), R(B) = Y,
N(B) = \{0\}$, $B$ is closed, and $F$ is finite-rank operator.

Let $\{\varphi_j\}_{1 \leq j \leq n}$ be a basis of $N(A)$ and
$\{\psi_j\}_{1 \leq j \leq n}$ be a basis of $N(A^\ast)$. It is known that
$$R(A)^\perp = N(A^\ast), \eqno{(2.1)}$$
where $R(A)^\perp$ is the set of linear functionals $\{\psi_j\}$
in $Y^\ast$ such that $(\psi_j, Au)=0$ $\forall u \in D(A)$, where
$(\psi_j, f)$ is the value of a linear functional $\psi_j \in Y^\ast$ on the
element $f \in Y$. Clearly, $\psi_j \in N(A^\ast), 1 \leq j \leq n$.

Define
$$Bu:=Au + \sum^n_{j=1} (h_j, u) \nu_j := A+F, \quad \nu_j \in Y,
  \eqno{(2.2)}$$
where $F$ is a finite-rank operator,
$\{\nu_j\}_{1 \leq j \leq n}$ is the set of elements of $Y$, biorthogonal
to the set $\{\psi_j\}_{1 \leq j \leq n}$,
$(\psi_j, \nu_m)= \delta_{jm} :=
 \begin{cases} 0, j \neq m \\ 1, j=m \end{cases}$,
and $\{h_j\}_{1 \leq j \leq n}$ is the set of elements of $X^\ast$,
biorthogonal to the set $\{\varphi_j\}_{1 \leq j \leq n}$,
$(h_j, \varphi_m) = \delta_{jm}$. Existence of the sets biorthogonal to
finitely many linearly independent elements of a Banach space follows from
the Hahn-Banach theorem. An arbitrary element $u \in X$ can be uniquely
represented as $u=u_1 + \sum^n_{j=1} c_j \varphi_j, c_j = \const$, and
$(h_j, u_1) = 0$.

Let us check that $N(B) = \{0\}$ and $R(B) = Y$. Assume $Bu=0$, that is
$Au + \sum^n_{j=1} (h_j, u) \nu_j =0$. Apply $\psi_m$ to this equation, use
$(\psi_m, Au) = 0$, and get
$$0 = \sum^n_{j=1} (\psi_m, \nu_j) (h_j, u) = \sum^n_{j=1} \delta_{mj}
(h_j, u) = (h_m, u), \quad 1 \leq m \leq n.$$
Therefore $Au=0$. So $u \in N(A),$ and $u=\sum^n_{j=1} c_j \varphi_j$,
$c_j = \const$.
Apply $h_m$ to this equation and use $(h_m, \varphi_j)= \delta_{mj}$ to get
$c_m=0$, $1 \leq m \leq n$. Thus $u=0$. We have proved that $N(B) = \{0\}$.

To prove $R(B) = Y$, take an arbitrary element $f \in Y$ and write
$f = f_1 + f_2$, where $f_1 = Au_1$ belongs to $R(A)$, and
$f_2 = \sum^n_{j=1} a_j \nu_j$, $a_j \equiv \const$. Note that
$$y = R(A) \dotplus L_n, \eqno{(2.3)}$$
where the sum is direct, $L_n$ is spanned by the elements
$\{\nu_j\}_{1\leq j \leq n}$, and $a_j = (\psi_j, f)$. Indeed,
$$(\psi_m, f) = (\psi_m, Au_1) + \sum^n_{j=1} a_j (\psi_m, \nu_j) = a_m,
  1 \leq m \leq n.$$
Given an arbitrary $f \in Y$, $f=Au_1 + \sum^n_{j=1} (\psi_j, f) \nu_j$,
define $u=u_1 + \sum^n_{j=1} (\psi_j, f) \varphi_j$. Then $Bu=f$. Indeed,
using (2.2) one has:
$$B\left[u_1 + \sum^n_{j=1}(\psi_j, f) \varphi_j \right] = Au_1 +
  \sum^n_{j=1} (h_j, u_1) \nu_j + \sum^n_{j=1} (h_j, \sum^n_{m=1}
  (\psi_m, f) \varphi_m) \nu_j = f. \eqno{(2.4)}$$

Here the relations $(h_j, \varphi_m) = \delta_{jm}$, 
and $(h_j, u_1) = 0$ are used.
We have proved the relation $R(B) = Y$.

2. Let us now assume that $A=B-F$, where $B:X\to Y$ is a linear closed densely
defined operator, $D(A) = D(B)$, $N(B) = \{0\}$, $R(B) = Y$, and $F$ is a
finite-rank operator. We wish to prove that (1.1) and (1.2) hold and $A$ is
closed.

Let us prove (1.1). Assume that $Au_n := f_n \to f$ and prove that
$f \in R(A)$.

One has $Bu_n-Fu_n \to f$. Since $N(B) = \{0\}$, $R(B) = Y$, and $B$ is
closed, $B^{-1}$ is bounded by Banach's theorem. Thus
$$u_n - B^{-1} Fu_n \to B^{-1} f. \eqno{(2.5)}$$
Since $F$ is a finite-rank operator,
$B^{-1}F$ is compact. Therefore, if $\sup_n \parallel u_n \parallel \leq c$,
then a subsequence, denoted $u_n$ again, can be found, such that
$B^{-1} Fu_n$ converges in the norm of $X$. Consequently, (2.5) implies
$u_n \to u$, $u-B^{-1} Fu = B^{-1} f$, so $u \in D(B)$ and
$Bu-Fu = f$.

To finish the proof, let us establish the estimate
$\sup_n \| u_n \| \leq c$. Assuming $\| u_{n_k} \| \to \infty$ and denoting
$n_k$ by $n$ and $B^{-1} F$ by $T$, define $v_n := \frac{u_n}{\|u_n\|}$,
$\| v_n \| = 1$. Then $v_n - T v_n \to 0$ as $n \to \infty$. One may assume
that $v_n$ is chosen in the direct complement of $N(I-T)$ in $X$. Arguing
as above, one selects a convergent in $X$ subsequence, denoted again by
$v_n$, $v_n \to v$, and gets $v-Tv=0$. Since $v$ belongs to the direct
complement of $N(I-T)$, it follows that $v=0$. On the other hand, since
$\| v \| = \lim_{n \to \infty} \| v_n \| = 1$, one gets a contradiction,
which proves the desired estimate $\sup_n \| u_n \| \leq c$.
Property (1.1) is proved. 

Let us prove that $A$ is closed. If $Au_n \to f$ and $u_n \to u$, then
$Bu_n - Fu_n \to f$, and the above argument shows that
$Bu - Fu = f$ so $Au = f$. Thus $A$ is closed.

Finally, let us prove (1.2).

Let $Au=0$, i.e. $Bu - Fu =0$. Applying a bounded linear injective operator
$B^{-1}$, one gets an equivalent equation
$$u - Tu = 0, \quad T:=B^{-1}F, \eqno{(2.6)}$$
with a finite-rank operator $T$. It is an elementary fact (see \cite{4}) that
$\dim N(I-T) := n < \infty$ if $T$ is a finite-rank operator.
Since $N(A) = N(I-T)$, one has $\dim N(A) = n < \infty$.

Now let $A^\ast v=0$. Then
$$B^\ast v-F^\ast v=0. \eqno{(2.7)}$$
Since $(B^\ast)^{-1} = (B^{-1})^\ast$ is a bounded and injective linear
operator, the elements $v$ are in one-to-one correspondence with the elements
$w:=B^\ast v$, and (2.7) is equivalent to
$$w-T^\ast w=0, \quad T^\ast = F^\ast (B^\ast)^{-1}, \eqno{(2.8)}$$
so that $T^\ast$ is the adjoint to operator $T:=B^{-1}F$.

Since $T$ is a finite rank operator, it is an elementary fact
(see \cite{4}) that
$\dim N(I-T^\ast) = \dim (I-T) = n < \infty$. Since $N(A^\ast) = N(I-T^\ast)$,
property (1.2) is proved.

Theorem 1 is proved. \qed

An immediate consequence of Theorem 1 is the Fredholm alternative (see
Theorem 1.1 in \cite{4}) for unbounded operators $A \in Fred(X,Y)$.

\end{document}